\documentclass[12pt]{amsart}
\usepackage{amsfonts}
\usepackage{latexsym,amssymb}
\usepackage{graphicx}
\usepackage[cp1250]{inputenc}
\usepackage[active]{srcltx}
\usepackage{amsaddr}

\newtheorem{theorem}{Theorem}
\newtheorem{corollary}[theorem]{Corollary}
\newtheorem{lemma}[theorem]{Lemma}

\theoremstyle{definition}
\newtheorem{remark}[theorem]{Remark}
\newtheorem{example}[theorem]{Example}

\newcommand{\Comp}{\mathbb{C}}

\newcommand{\Zzz}{\mathbb{Z}}

\newcommand{\cB}{\mathcal{B}}
\newcommand{\cH}{\mathcal{H}}
\newcommand{\cK}{\mathcal{K}}
\newcommand{\cM}{\mathcal{M}}
\newcommand{\cR}{\mathcal{R}}

\newcommand{\bC}{\mathbb{C}}
\newcommand{\bN}{\mathbb{N}}

\newcommand{\li}{l^\infty}

\newcommand{\matp}[1]{\begin{bmatrix} #1 \end{bmatrix}}

\usepackage{color}

   \author{Ahmed Bachir}
   \address{
Department of Mathematics\\
King Khalid University\\
Abha\\
Saudi Arabia\\}
   \email{bachir\_ahmed@hotmail.com}

	\author{Patryk Pagacz}
  \address{(Corresponding author)\\
   Faculty of Mathematics and Computer Science\\
   Jagiellonian University\\
   \L ojasiewicza~6\\
   30-348 Krak\'ow\\
   Poland}
   \email{patryk.pagacz@gmail.com}

\begin{document}

   \title[  ]
   {The Berberian's transform and an asymmetric Putnam-Fuglede theorem}
	
	

\date{}%
\begin{abstract}
In this paper we present how to apply a Berberian's technique to asymmetric Putnam-Fuglede theorems. In particular, we proved that if $A,B\in \cB(\cH)$ belong to the union of classes of $*$-paranormal operators, $p$-hyponormal operators, dominant operators and operators of class $\mathcal{Y}$ and $AX=XB^*$ for some $X\in \cB(\cH)$, then $A^{*}X=XB$. 

Moreover, we gave a new counterexample for an asymmetric Putnam-Fuglede theorem for paranormal operators.
\end{abstract}

\maketitle

\section{Introduction}

The familiar Putnam-Fuglede's theorem asserts that if $A\in \cB(\cH)$ and $B \in \cB(\cH)$ are normal operators and $AX=XB$ for some $X\in \cB(\cH)$, then $A^{*}X=XB^{*}$ (see \cite{PUTN}). A simple example with $A=B=X$ being a unilateral shift shows that this implication cannot be extended to the class of subnormal operators. Thus, let us overwrite the Putnam-Fuglede's theorem in an asymmetric form: 

\begin{theorem}[Putnam-Fuglede's theorem]
Let $A,B \in \cB(\cH)$ be normal and $AX=XB^*$ for some $X\in \cB(\cH)$, then $A^{*}X=XB$.
\end{theorem}
 Many authors extended this theorem for several classes consisting subnormal operators.

At the beginning, Stampfli and Wadhwa proved that the asymmetric Putnam-Fuglede's theorem holds for dominant and normal operators, i.e. they assumed that $A$ is dominant and $B$ is normal (see \cite{STAM}). Later they extend the result to dominant and hyponormal operators (see \cite{STAM2}).

Independently, Futura in \cite{F2} shows that  the asymmetric Putnam-Fuglede's theorem holds for both subnormal operators, i.e. he assumed that $A$ and $B$ are subnormal. The analogy theorem for unbounded operators was proved by Stochel (see \cite{S}).

Next, a three important papers made a further extensions.
Namely, Gupta in \cite{Gupta} shows the asymmetric Putnam-Fuglede's theorem for both $k$-quasihyponormal injections.
Moore, Rogers and Trent in \cite{MRT} shows the theorem for both M-hyponormal operators. Yoshino in \cite{Y} (and Duggal in \cite{D}) proved it when $A$ is M-hyponormal and $B$ is dominant operator.

A few years later, Uchiyama and Takanashi in \cite{UchiyamaT} proved the asymmetric Putnam-Fuglede's theorem in three cases: when $A$ is $p$-hyponormal or $log$-hyponormal, $B$ is dominant, when both $A$ and $B$ are $p$-hyponormals, and when  both $A$ and $B$ are $log$-hyponormals.

Recently many authors (f. e. see \cite{FX,Kim,R,Rashid,UchiyamaY}) try to extend the Putnam-Fuglede theorem for superclasses of normal operators. Other direction of developing the Putnam-Fuglede theorem, started by Berberian(\cite{B_H-S}), is to put an additional assumption onto interwinding operator f.e. to assume that $X$ belongs to the Hilbert-Schmidt class, i.e. $X^*X$ has a finite trace (f.e. see \cite{F3,MU}).

The organization of the paper is as follows, in Section 2, we give all necessary definitions and we illustrate the relations between the main superclasses of normal operators. Moreover, we briefly recall Berberian's construction.

In Section 3, we highlight the operators for which the only invariant subspace which restrict the operator to a normal one is reducing \eqref{eq:main}. We prove that if Berberian's extensions of $A$ and $B$ satisfy this property then $A$ and $B$ satisfy the asymmetric Putnam-Fuglede theorem. This theorem can be easily used for some superclasses of normal operators.

In Section 4, we use our tool theorem to show the asymmetric Putnam-Fuglede's theorem for some certain superclasses of normal operators. This extends the recent results form \cite{DJK}, \cite{FX} and \cite{P}.

Since all powers of hyponormal operator are paranormal (see \cite{F}), but not necessary hyponormal (see \cite{H,R}), many authors gave their attention to paranormal operators. Moreover, most of operators considered in this article are paranormal we also bring our attention in Section 5 to paranormal operators. 
In \cite{R} Radjabalipour showed that for each $A$ and $B$ such that $A,B$ are the same powers of some hyponormal operators the asymmetric Putnam-Fuglede's theorem holds true. The same is not true for paranormal operators (see \cite{P,Uchiyama}). 

In \cite{Uchiyama} Uchiyama showed an example of paranormal operator with non-reducing eigenspace. In Section 5 we gave a simpler example which shows that an asymmetric Putnam-Fuglede's theorem for paranormal and unitary operator does not hold, even if one assume that the interwinding operator $X$ belongs to Hilbert-Schmidt class. This may surprise, since Mecheri and Uchiyama in \cite{MU} showed that if $A$ is a class A operator, $B$ is an injective class A operator and $X$ is a Hilbert-Schmidt operator then the Putnam-Fuglede's theorem holds true.

\section{Preliminaries}
Throughout this paper,  $\cH $ denotes an infinite dimensional complex Hilbert space with inner product $\langle \cdot, - \rangle$  and  $\cB(\cH)$  denotes the algebra of all bounded linear operators acting on $\cH$. Spectrum, point spectrum, residual spectrum, continuous spectrum and approximate spectrum of an operator $T$ will be denote by $\sigma(T)$, $\sigma_p(T)$, $\sigma_r(T)$, $\sigma_c(T)$, $\sigma_{ap}(T)$, respectively. The kernel and the range of an operator $T$ will be denote by $\ker T$ and $\cR(T)$. The subspace $\cK \subset \cH$ is \emph{invariant} for $T\in\cB(\cH)$ if $T(\cK)\subset \cK$ and \emph{reducing} for $T$ if $T(\cK) \subset \cK$ and $T^*(\cK) \subset \cK$.

For any operator $T\in \cB(\cH)$, set, as usual $|T|=(T^{*}T)^{1/2}$ and $[T^{*},T]=T^{*}T-TT^{*}$ (the self-commutator of $T$), and  consider the following standard definitions:
\begin{itemize}\label{alldefs}
  \item $T$ is \emph{normal} if $T^*T=TT^*$,
\item $T$ is \emph{hyponormal} if $|T^{*}|^{2}\leq |T|^{2}$ (i.e. if $[T^{*},T]$ is nonnegative or equivalently, if $\|T^{*}h\|\leq \|Th\|$ for every $h\in \cH$),
\item $T$ is \emph{M-hyponormal} if there exists $M>0$ such that $\|(T^*-\overline{\lambda})h\|\leq M \|(T-\lambda) h\|$, for any $h\in\cH$ and $\lambda\in\Comp$, see \cite{Wadhwa},
\item $T$ is \emph{dominant} if for any $\lambda\in\Comp$ there exists $M_{\lambda}>0$ such that $\|(T^*-\overline{\lambda})h\|\leq M_{\lambda} \|(T-\lambda) h\|$, for any $h\in\cH$, see \cite{STAM},
\item $T$ is \emph{$p$-hyponormal} if $(T^*T)^p\geq (TT^*)^p$, see \cite{Alu0},
\item $T$ is \emph{log-hyponormal} if $T$ is invertible and $\log(T^*T)\geq \log(TT^*)$, see \cite{Tana},
\item $T$ is \emph{a class A} operator if $|T^2|\geq |T|^2$, see \cite{F4},
\item $T$ is \emph{of class $\mathcal{Y}$} if for any $\alpha\geq 0$ there exists $k_\alpha$ such that $|T^*T-TT^*|^{\alpha}\leq k_\alpha^2(T-\lambda)^*(T-\lambda)$ for all $\lambda \in\Comp$, see \cite{UchiyamaY},

\item $T$ is \emph{paranormal} if $T^{*2}T^2-2\lambda T^*T+\lambda^2\geq 0$, for all $\lambda>0$ (or equivalently, if $\|Th\|^2\leq \|T^2h\|$ for all $h\in \cH$), see \cite{F},
\item $T$ is \emph{$*$-paranormal} if $T^{*2}T^2-2\lambda TT^*+\lambda^2\geq 0$, for all $\lambda>0$ (or equivalently, if $\|T^{*}h\|^2\leq \|T^2h\|$ for every $h\in \cH$), see \cite{PhD}.
\end{itemize}

The inclusion relations between the above-mentioned operators classes are shown
below.
$$hyponormal\subset M-hyponormal \subset dominant$$
$$hyponormal\subset p-hyponormal \subset A \ class \subset paranormal$$
$$invertible \ p-hyponormal \subset \log-hyponormal \subset A \ class $$

Finally, let us remain the Berberian's technique.

Let $\li(\cH)$ denote the space of all bounded sequences of elements of $\cH$. Let us fix $\phi$ a Banach limit. Let $c_0(\cH)$ denote the space of all null sequences of $\cH$, i.e. $c_0(\cH):=\{\{x_n\}_n \subset \cH : \phi(\{\|x_n\|^2\}_n)=0\}$. Endowed with the canonical norm, the quotient space $\cK = \li(\cH )/c_0(\cH)$ can be made into a Hilbert space (see \cite{B}). The transform $\cH\ni x \mapsto \{x\}_{n\in \bN} + c_0(\cH)\in \cK$ is a natural isometric embedding. By \cite{B} there exists an isometric $*$-isomorphism $\cB(\cH)\ni T \to T^\circ \in \cB(\cK)$ preserving order such that $\sigma(T)= \sigma(T^\circ)$ and $\sigma_{ap}(T)=\sigma_{ap}(T^\circ)=\sigma_{p}(T^\circ)$.

Moreover, by the uniqueness of the square root of a positive operator, one can deduce that $(A^\frac12)^\circ=(A^\circ)^\frac12$ for any positive $A\in\cB(\cH)$.
Thus it is clear that $|T^\circ|=|T|^\circ$ for any $T\in\cB(\cH)$.

\section{Main Theorem}

Let $T\in\cB(\cH)$ be an operator which \emph{satisfies the Putnam-Fuglede theorem} i.e. for all operators $X,N\in\cB(\cH)$ such that $N$ is normal and $TX=XN$, it holds that $T^*X=XN^*$. Then $T$ shares the following property
\begin{equation}\label{eq:main}\tag{$\spadesuit$}
\textnormal{any invariant subspace }\cM \subset \cH\textnormal{ of }T \textnormal{ such that }T|_{\cM}\textnormal{ is normal,}
\end{equation}
$$\textnormal{ is reducing for } T.$$
Indeed, if a subspace $\cM \subset \cH$ is invariant for $T$ and $T|_{\cM}$ is normal, then $TX=XN$, where $X=P_{\cM}$ and $N=T|_{\cM}\oplus Id_{\cM^\bot}$. Thus $T^*X=XN^*$. Hence $(T|_{\cM})^*=T^*|_{\cM}$ and so $\cM$ is reducing for $T$.

On the other hand, not all operators which share the property \eqref{eq:main} satisfy Putnam-Fuglede theorem.

\begin{example}
Let consider a Hilbert space $l^2(\mathbb{Z})$ with a canonical basis $\{e_n\}_{n\in\Zzz}$. Let $S,N$ be bounded operators define on the basis as
follows
$Se_n=\begin{cases}
   e_{n+1}, & \mbox{if } n\not=-1 \\
   2e_0, & \mbox{if } n=-1
 \end{cases}$; $Ne_n = e_{n+1}$ for all $n \in \Zzz$. The operator $N$ is normal. Moreover, let take $X = \frac12 Id_{l^2(\mathbb{Z}_-)}\oplus Id_{l^2(\mathbb{Z}_+)}$, where
$\Zzz_+ := \{ 0,1, 2,\dots\}$, $\Zzz_+ := \{-1,-2,-3,\dots \}$.

We have $SXe_n = XNe_n = \begin{cases}
                           \frac12 e_n, & \mbox{if } n<-1 \\
                           e_n, & \mbox{if } n \geq -1.
                         \end{cases}$,
thus $SX = XN$
But $S^*Xe_0=2 e_{-1} \not= \frac12 e_{-1} = XN^*e_0$, hence $S^*X \not= XN^*$.

Thus $S$ does not satisfy Putnam-Fuglede theorem and $S$ satisfy \eqref{eq:main}, since it cannot be restricted to a normal operator.
\end{example}

Let us remark a simple fact about property \eqref{eq:main}.

\begin{remark}\label{simpleRemark}
Let $T\in\cB(\cH)$.
\begin{itemize}
  \item If $T$ is normal, then it satisfies \eqref{eq:main}.
  \item If $T$ satisfies \eqref{eq:main}, then for any unitary operator $U\in\cB(\cH)$ the operator $UTU^*$ satisfies \eqref{eq:main}.
  \item If $T$ satisfies \eqref{eq:main}, then for any invariant subspace $\cM\subset \cH$ the operator $T|_{\cM}$ satisfies \eqref{eq:main}.
  \item If $T^\circ$ satisfies \eqref{eq:main}, then $T$ satisfies \eqref{eq:main}.
\end{itemize}
\end{remark}

Before we prove our main theorem let us observe how the property \eqref{eq:main} impact the spectrum of the adjoint.

\begin{lemma}\label{cztery}
Let $T\in\cB(\cH)$ be an operator which satisfy \eqref{eq:main}, then the residual spectrum of $T^*$ is empty. In particular, the spectrum  of $T^*$ is equal to approximate spectrum of $T^*$.
\end{lemma}
\begin{proof}
The space $\cH$ is a direct sum of the space spanned by eigenvectors of $T$ and its complement. Thus by \eqref{eq:main} $T$ is a direct sum of scalar operators and an operator $T'$ with empty point spectrum. Then since $ker (\lambda - T')=\{0\}$, we have $\overline{ \cR(\overline{\lambda} - T'^*)} = \cH$ for all $\lambda\in \bC$.
Hence $\sigma_r(T'^*)=\emptyset$ and $\sigma_r(T^*)=\emptyset$.
Moreover, $\sigma(T^*) = \sigma_p(T^*)\cup \sigma_c(T^*) \subset \sigma_{ap}(T^*)$.
\end{proof}

Our main theorem shows the importance of property \eqref{eq:main}.

\begin{theorem}\label{mainTheorem}
Let $A,B \in \cB(\cH)$ be such that $A^\circ,B^\circ$ satisfy \eqref{eq:main}.
If $AX = XB^*$, for some $X\in \cB(\cH)$, then $A^*X = XB$.
\end{theorem}

\begin{proof}
Let us fixed $X\in \cB(\cH)$ and consider its polar decomposition $X=U|X|$, with a unitary operator $U: \overline{\cR(|X|)} \to \overline{\cR(X)}$. Then the equation $AX = XB^*$ is equivalent to $\widetilde{A}|X| = |X|B^*$, where $\widetilde{A}:=U^{-1}AU \oplus 0_{ker |X|}$. Thus by Remark \ref{simpleRemark} it is enough to prove the implication for a positive operator $X$.

Thus let operators $A^\circ,B^\circ$ satisfy \eqref{eq:main} and positive operator $X$ be such that $AX = XB^*$. Hence the subspace $\overline{\cR(X)}$ is invariant for $A$. Since the subspace $ker X$ is invariant for $B^*$, then $\overline{\cR(X)}$ is also invariant for $B$. As a consequence we have the following matrices representations with respect to the decomposition $\cH= \overline{\cR(X)} \oplus ker X$.
$$
X=\matp{K & 0 \\ 0 & 0},\quad A=\matp{A_{11} & A_{12} \\ 0 & A_{22}}, \quad B=\matp{B_{11} & B_{12} \\ 0 & B_{22}}.
$$

The equation $AX = XB^*$ implies $A_{11}K = KB_{11}^*$. By Lemma \ref{cztery} we have $\sigma(B_{11}^*) = \sigma_{ap}(B_{11}^*)$.

The Berberian's extensions $A_{11}^\circ, B_{11}^\circ, K^\circ$ of the operators $A_{11}, B_{11}, K$ satisfy the equation
\begin{equation}\label{eq}
A^\circ_{11}K^\circ = K^\circ (B_{11}^*)^\circ
\end{equation}
and $\sigma((B_{11}^*)^\circ)=\sigma(B_{11}^*) =\sigma_{ap}(B_{11}^*)=\sigma_{p}((B_{11}^*)^\circ)$.
The equation \eqref{eq} is equivalent to $$(\lambda - A_{11}^\circ) K^\circ  = K^\circ (\lambda - (B^*_{11})^\circ),$$ for $\lambda \in \bC$. Thus if $\lambda \in \sigma_r(A_{11}^\circ)$, then $\lambda \in \sigma_r((B^*_{11})^\circ)$. But  by Lemma \ref{cztery} we get  $\sigma_r((B_{11}^*)^\circ)=\emptyset$. As a consequence $\sigma_r(A_{11}^\circ)=\emptyset$. So $\sigma(A_{11}^\circ)=\sigma_{ap}(A_{11}^\circ)=\sigma_p(A_{11}^\circ)$. Moreover, the operator $A_{11}^\circ$ satisfies \eqref{eq:main} thus it is a direct sum of scalar operators.
Normality of $A_{11}^\circ$ shows that $A_{11}$ is normal. Hence by \eqref{eq:main} for $A$ (see Remark \ref{simpleRemark}) we get $A_{12} =0$.

The equation \eqref{eq} is equivalent to $K^\circ (A_{11}^*)^\circ = B_{11}^\circ K^\circ $. Thus we can repeat the above argument and show that the operator $B_{11}$ is normal and $B_{12} =0$.

Finally, to show that $A^*X=XB$ it is enough to show that $A_{11}^*K=KB_{11}$, but it is a consequence of the classical Putnam-Fuglede theorem.
\end{proof}

We can resume this section in the following

\begin{corollary}
Let $T\in\cB(\cH)$. Then
\begin{itemize}
  \item $T^\circ$ satisfies \eqref{eq:main}, if $T$ is normal,
  \item $T$ satisfies Putnam-Fuglede, if $T^\circ$ satisfies \eqref{eq:main},
  \item $T$ satisfies \eqref{eq:main}, if $T$ satisfies Putnam-Fuglede.
\end{itemize}
\end{corollary}

\section{Putnam-Fuglede theorem for some classes of operators }

In the section we will show the property \eqref{eq:main} for some classes of operators. Thus by Theorem \ref{mainTheorem} we prove a Putnam-Fuglede asymmetric theorem.

\begin{lemma}\label{pik_*para}
Any $*$-paranormal operator satisfies \eqref{eq:main}.
\end{lemma}
\begin{proof}
  Let $T\in\cB(\cH)$ be a $*$-paranormal operator and $\cM$ be a subspace invariant for $T$ such that  $T|_{\cM}$ is normal.
  In other words,  $T=\begin{pmatrix} N & A \\ 0 & B \end{pmatrix}$ on $\cH=\cM \oplus \cM^\perp$, where $N:=T|_{\cM}$ is a normal operator on $\cM$.
	
The $\ast$-paranormality of $T$ implies that
$$T^{2*}T^2 -2k TT^*+k^2I \geq 0, \quad \textnormal{for all } k>0.$$
Hence
$$|N|^4-2k (|N|^2 +AA^*)+k^2I_{\cM}\geq 0, \quad \textnormal{for all } k>0.$$
Thus
$$
2k AA^* \leq (|N|^2-kI_{\cM})^2, \quad \textnormal{for all } k>0,
$$
and by Douglas's theorem
$$\cR(A) \subset \cR(|N|^2 -kI_{\cM}), \quad \textnormal{for all } k>0.$$
Since $|N|^2-\lambda I_{\cM}$ is invertible for all $\lambda \in \mathbb{C}\setminus [0, \infty)$, we have
$$\cR(A) \subset \bigcap_{\lambda\not= 0}\cR(|N|^2-\lambda  I_{\cM}).$$
Furthermore, by Theorem 1 in \cite{Putnam} one can get
$$ \bigcap_{\lambda\not= 0}\cR(|N|^2-\lambda  I_{\cM}) =\cR(E(\{0\}))= \ker |N|^2=\ker N,$$
where $E$ is a spectral measure of the operator $|N|^2$.

On the other hand,
$$\|A^*x\|^2=\|N^*x\|^2+\|A^*x\|^2=\|T^*x\|^2\leq \|T^2x\|\|x\|=\|N^2x\|\|x\|=0,$$ for any $x\in \ker N \subset \cM$. Hence $\ker N^*=\ker N \subset \ker A^*$. Thus $\overline{\cR(A)}\subset \overline{\cR(N)}$ and $\overline{\cR(A)}\subset \ker N$.
So $A=0$. Therefore, $\cM$ is a subspace  reducing for $T$.
\end{proof}

The analogy statement to Lemma \ref{pik_*para} is true for other superclasses of normal operators (see \cite{Rad},\cite{UchiyamaT},\cite{UchiyamaY}).
\begin{lemma}\label{othercirc}
The property \eqref{eq:main} is satisfied for
\begin{itemize}
  \item $p$-hyponormal operators,
  \item dominant operators,
  \item operators of class $\mathcal{Y}$.
\end{itemize}
\end{lemma}

By the above lemmas Theorem \ref{mainTheorem} can be apply for several classes of operators.

\begin{theorem}
Let $A,B\in \cB(\cH)$. Let $A$ as well as $B$ be 
\begin{itemize}
\item $*$-paranormal operator or
\item $p$-hyponormal operator, where $1>p>0$, or
\item dominant operator or
\item operator of class $\mathcal{Y}$.
\end{itemize}
If $AX = XB^*$, for some $X\in \cB(\cH)$, then $A^*X = XB$.
\end{theorem}

\begin{proof}
Let us recall that $(C^\circ)^\frac12=(C^\frac12)^\circ,$ for any positive operator $C\in\cB(\cH)$.
As a consequence $(C^\circ)^p=(C^p)^\circ$, for any $1>p>0$ with finite binary representation. 
But, since the function $p\mapsto \langle A^px,x \rangle$ is convex, the above equality holds for any $1>p>0$.
Thus, by Definition \ref{alldefs} and properties of the Berberian's extension, operators $A^\circ,B^\circ$ also belong to the one of the above-mentioned classes. As a consequence of lemmas \ref{pik_*para},\ref{othercirc}, the assumptions of Theorem \ref{mainTheorem} are satisfy.
\end{proof}

\begin{remark}
One can find the above theorem for *-paranormal operators in \cite{Rashid} (Section 5). But unfortunately that part of \cite{Rashid} followed the unpublished, incorrect version of a recent article (\cite{PBarxiv}).
\end{remark}

\section{Putnam-Fuglede theorem for paranormal operators modulo Hilbert-Schmidt operators}

In \cite{P} we presented an example of paranormal operator $S$, a unitary operator $U$ and orthogonal projection $P$ such that $SP=PU$, but $S^*P\not=PU^*$. Thus we show that an asymmetric Putnam-Fuglede theorem for paranormal operators does not hold. But ealier, Uchiyama and Tahanashi in \cite{UchiyamaT} gave an example of a Hilbert space with the basis $\{e_{n}\}_{n=-\infty}^\infty \cup \{f\}$ and $A$ class operator $T$ such that $\ker T=\Comp (f-e_{-1})$ and $T^*(f-e_{-1})=-e_{-2}$. This shows that for $P=P|_{\ker T}$ and $N=P_{\Comp e_{-17}}$ we have $TP=0=PN$, but $T^*P\not=0=PN^*$. Later in \cite{MU} Uchiyama and Mecheri showed that $N$ cannot be taken as an injective operator. In fact, they proved that an asymmetric Putnam-Fuglede theorem holds modulo Hilbert-Schmidt operators for $A$ class operators, if we assume that one of them is injective.

Here we improve the example from \cite{P} to the case where $P$ is one-dimension projection (a Hilbert-Schmidt operator).
Therefore, the result from \cite{UchiyamaT} for $A$ class operators cannot be extend to paranormal operators.

\begin{example}
Let 
$$T:l^2\ni (x_0,x_1,x_2,\dots)\mapsto (x_0+x_1,x_1,x_1,\sqrt{8}x_2,\sqrt{8}x_3,\dots) \in l^2.$$
This operator can be express as follows
$$T=\matp{1 & 1 & 0 & 0 & 0 & \dots \\ 0 & 1 & 0 & 0 & 0 & \dots \\ 0 & 1 & 0 & 0 & 0 & \dots \\ 0 & 0 & \sqrt{8} & 0 & 0 & \dots \\ 0 & 0 & 0 & \sqrt{8} & 0 & \dots \\ \vdots & \vdots & \vdots & \vdots & \vdots & \ddots } \in \cB(l^2).$$
The operator $T$ is paranormal.

Let $h\in l^2$, then $h=\alpha e_0 + \beta e_1 + \sum\limits_{n=2}^\infty \gamma_n e_n$, where $\{e_n\}_{n\in\bN}$ is an orthogonal basis of $l^2$. Since three sets $\{\alpha e_0 + \beta e_1, e_2, e_3,\dots\}$, $\{T(\alpha e_0 + \beta e_1), Te_2, Te_3,\dots\}$ and $\{T^2(\alpha e_0 + \beta e_1), T^2e_2, T^2e_3,\dots\}$ consist of orthogonal vectors, we have $$\|T^2h\|^2-2\lambda\|Th\|^2+\lambda^2\|h\|^2= \|T^2(\alpha e_0 + \beta e_1)\|^2-2\lambda\|T(\alpha e_0 + \beta e_1)\|^2+\lambda^2\|\alpha e_0 + \beta e_1\|^2+$$ $$+\sum\limits_{n=2}^\infty\|T^2\gamma_ne_n\|^2-2\lambda\|T\gamma_ne_n\|^2+\lambda^2\|\gamma_ne_n\|^2.$$
Now let us observe that
$$\|T^2\gamma_ne_n\|^2-2\lambda\|T\gamma_ne_n\|^2+\lambda^2\|\gamma_ne_n\|^2 = \|\gamma_ne_n\|^2(\lambda-8)^2\geq0.$$
Moreover,
$$\|T^2(\alpha e_0 + \beta e_1)\|^2-2\lambda\|T(\alpha e_0 + \beta e_1)\|^2+\lambda^2\|\alpha e_0 + \beta e_1\|^2=$$
$$=|\alpha e_0+\beta+\beta|^2+|\beta|^2+8|\beta|^2-2\lambda(|\alpha e_0+\beta|^2+2|\beta|^2)+\lambda(|\alpha|^2+|\beta|^2)=$$
$$=|(1-\lambda)\alpha+2\beta|^2+(3-t)^2|\beta|^2\geq 0.$$
Thus for every positive number $\lambda$ and $h\in\cH$ we get $\|T^2h\|^2-2\lambda\|Th\|^2+\lambda^2\|h\|^2 \geq 0$. Hence $T$ is paranormal.
If we take $P$ an projection onto one-dimension space $\bC e_0$ and $U=Id$, then we get $TP=PU$, but $T^*P e_0= e_0 + e_1 \not= e_0 = PU^* e_0$.
\end{example}


\begin{thebibliography}{99}

\bibitem{Alu0} A. Aluthge, \textit{On p-hyponormal operators for $0 < p < 1$}, Integral Equations Operator Theory, 13(1990), 307--315.

\bibitem{B} S.K. Berberian, \textit{Approximate proper vectors}, Proc. Amer. Math. Soc., 13(1962), 111--114.

\bibitem{B_H-S} S.K. Berberian, \textit{Extensions of a theorem of Fuglede and Putnam}, Proc. Amer. Math. Soc., 71(1978), 113--114.

\bibitem{PBarxiv} A. Bachir, P. Pagacz \textit{An asymmetric Putnam-Fuglede theorem for $*$-paranormal operators}, arXiv:1405.4844v1 [math. FA], 19 May 2014.

\bibitem{D} B.P. Duggal, \textit{On dominant operators}, Arch. Math., 46(1986), 353--359.

\bibitem{DJK} B.P. Duggal, I.H. Jeon, I.H. Kim, \textit{On $*$-paranormal contractions and properties for $*$-class A operators}, Linear Algebra Appl.,  436 (2012), 954--962.

\bibitem{F} T. Furuta, \textit{On the class of Paranormal Operators}, Proc. Japan Acad., 43(1967), 567--649.

\bibitem{F2} T. Furuta, \textit{On relaxation of normality inthe Fuglede-Putnam Theorem}, Proc. Amer. Math. Soc., 77(1979), 324--328.

\bibitem{F3} T. Furuta, \textit{An extension of the Fuglede-Putnam Theorem to subnormal operators using a Hilbert-Schmidt norm inequality}, Proc. Amer. Math. Soc., 81(1981), 240--242.

\bibitem{F4} T. Furuta, M. Ito and T. Yamazaki, \textit{A subclass of paranormal operators including class of log-hyponormal and several related classes}, Scientiae Mathematicae 1(1998), 389--403.

\bibitem{FX} F. Gao, X. Li, \textit{On $*$-class A contractions}, J. Inequal. Appl.,  2013, 2013--239.

\bibitem{Gupta} B.C. Gupta, \textit{An extension of the Fuglede-Putnam theorem and normality of operators}, Indian J. Pure Appl. Math., 14(1983), 1343--1347.


\bibitem{H} P.R. Halmos, \textit{Hilbert Space Problem Book}, Van Nostrand. The University Series in Higher Mathematics (1966).

\bibitem{Kim} I.H.Kim, \textit{The Fuglede-Putnam theorem for $(p,k)$-quasihyponormal operators}, J.Ineq. Appl, (2006), article ID 47481, pp. 1--7.

\bibitem{MRT} R.L. Moore, D.D. Rogers, T.T. Trent, \textit{A note on intertwining M -hyponormal operators}, Proc. Amer. Math. Soc., 83(1981), 514--516.

\bibitem{MU} S. Mecheri, A. Uchiyama, \textit{An extension of the Fuglede-Putnam's theorem to class A operators}, Math. Ineq. Appl., 13(2010), 57--61.

\bibitem{STAM} J.G. Stamfli, B. Wadhwa, \textit{An asymmetric Putnam-Fuglede theorem for dominant operators}, Indian Univ. Math. J., 25, no. 4(1976), 359--365.

\bibitem{STAM2} J.G.Stampfli, B.L. Wadhwa. \textit{On dominant operators}, Monatsshefte f$\ddot{u}$r Math., 84(1977), 143--153.

\bibitem{S} J. Stochel, \textit{An asymmetric Putnam-Fuglede theorem for unbounded operators}, Proc. Amer. Math. Soc.,  129(2001), 2261--2271.

\bibitem{P} P. Pagacz, \textit{The Putnam-Fuglede property for paranormal and *-paranormal operators}, Opuscula Math. 33, no. 3 (2013), 565--574.

\bibitem{PhD} S.M. Patel, \textit{Contributions to the study of spectraliod operators}, PhD thesis,
Delhi University, 1974.

\bibitem{PUTN} C.R. Putnam, \textit{On the normal operators on Hilbert space}, Amer. J. Math., 73(1951), 357--362.

\bibitem{Putnam} C.R. Putnam, \textit{Range of normal and subnormal operators}, Michigan Math. J., 18(1971), 33--36.

\bibitem{Rad} M. Radjabalipour, \textit{Majorization and normality of operators}, Mathematische Zeitschrift, Proc. Amer. Math. Soc., 62(1977), 105--110.

\bibitem{R} M. Radjabalipour, \textit{An extension of Putnam-Fuglede Theorem for Hyponormal Operators}, Mathematische Zeitschrift, 194(1987), 117--120.

\bibitem{Rashid} M. Rashid, \textit{On n$*$-paranormal operators}, Commun. Korean Math. Soc. 31(2016), 549--565.

\bibitem{Tana} K. Tanahashi, \textit{On log-hyponormal operators}, Integral Equations Operator Theory, 34(1999), 364--372.

\bibitem{UchiyamaT} K. Tahanashi and A. Uchiyama, \textit{Fuglede-Putnam's theorem for p-hyponormal or log-hyponormal operatos}, Glasgow Math. J., 44(2002), 397--410.


\bibitem{UchiyamaY} A. Uchiyama and T. Yoshino,\textit{ On the class $\mathcal{Y}$ operators}, Nihonkai Math. J., 8(1997), 179--194.

\bibitem{Uchiyama} A. Uchiyama,\textit{ An Example of Non-Reducing Eigenspace of Paranormal Operator}, Nihonkai Math. J., 14(2003), 121--123.

\bibitem{Wadhwa} B. L. Wadhwa, \textit{M-hyponormal operators}, Duke Math. Jour., 41(1974) 655--660.


\bibitem{Y} T. Yoshino, \textit{Remark on the generalized Putnam--Fuglede theorem}, Proc. Amer. Math. Soc., 95(1985), 571--572.





\end{thebibliography}
\end{document}